\documentstyle[12pt]{article}
\def\CC{{\rm\kern.24em\vrule width.02em height1.4ex depth-.05ex\kern-.26em
C}} 
\def\QQ{{\rm\kern.24em\vrule width.02em height1.4ex depth-.05ex\kern-.26em
Q}}
\def\PP{{\rm\kern.24em\vrule width.02em height1.4ex depth-.05ex\kern-.26em
P}}
\def\Rr{{\rm I\kern-.2em R}}
\def\ZZ{{\rm\kern.26em\vrule width.02em height0.5ex
depth0ex\kern.04em\vrule width.02em height1.47ex depth-1ex\kern-.34em Z}}
\def\BB{{\rm\kern.24em\vrule width.02em height1.4ex depth-.05ex\kern-.26em
B}}
\def\RR{\hspace{.065in}\rm{\vrule width.02em height1.55ex
depth-.07ex\kern-.3165em R}}
\def\Ibb#1{{\rm I\kern-.23em#1}}

\newcommand{\bea}{\begin{eqnarray*}}
\newcommand{\eea}{\end{eqnarray*}}
\def\CQFD{\hfill \vrule width 7pt height 7pt depth 1pt}
\newtheorem{lem}{LEMMA}[section]
\newtheorem{theo}[lem]{THEOREM}

\newtheorem{coro}[lem]{COROLLARY}

\newtheorem{remark}[lem]{Remark}

\begin{document}

In neutral/res*/fac/fornaess/qu*/papersi*/qu*july99

\medskip

\centerline{\Large \bf Infinite dimensional Complex Dynamics}

\bigskip

\centerline{Quasiconjugacies, Localization and Quantum Chaos}

\bigskip

\centerline{by}

\bigskip

\centerline{J. E. Forn\ae ss}

\bigskip

Version of
\today

\bigskip

\tableofcontents

\section{Introduction}

\indent The theory of complex dynamics started in the 1870's with works
by Schroder ([S]) in one complex dimension. Complex dynamics in
$\CC^2$ began with Fatou ([Fa]) in the second decade of this century. 
There has been
a lot of work in finite-dimensional complex dynamics, especially
in the last two decades, see for example ([F]) and references therein.

\medskip

In this paper we will make a first step towards a theory
of complex dynamics in infinite dimension. The main initial
problem is to define the basic concepts such as Chaotic
behaviour etc. For this we will
center our discussion around a concrete example. We also need to
find an example which is sufficiently elementary,
mathematically simple and easily modified to create distinct
dynamical behaviour, so that it can be used for testing 
concepts. Hence we will
find an infinite dimensional complex manifold $M$,
a holomorphic map $F:M \rightarrow M$ and investigate 
chaotic features, or lack thereof. One natural way to measure the dynamics
is to compare with a finite-dimensional approximation,
$f:N \rightarrow N$ since finite-dimensional systems have
been extensively studied. Another advantage to this comparison
is that infinite dimensional manifolds are not locally compact
while comparisons with the finite dimensional space
provides a useful substitute.

\medskip

The maps $f$ and $F$ can be compared by defining a map
$\Lambda : M \rightarrow N$ and comparing $f^n \circ \Lambda$
and $\Lambda \circ F^n.$ We call such $\Lambda$ quasiconjugacies.

\medskip

We choose our example from the field of 
Quantum Chaos. This is very
natural because solutions of the Schrodinger equation
are holomorphic maps $F:M \rightarrow M$ on
infinite-dimensional complex projective Hilbert space,
and can be compared with finite-dimensional maps
$f: N \rightarrow N$ describing the classical case.
Quantum Chaos is also a reasonable topic since it is a field without
well defined dynamical concepts and often without mathematical rigour.

\medskip

To be as specific as possible, we will focus on one of the most
basic problems of dynamics, namely
whether orbits of $F$ and $f$ are bounded. The phenomenon
that some orbits of $f$ are unbounded while all orbits
of $F$ are bounded is called localization. Our example is of
this type, see Theorem 4.5 and Corollary 4.8. 

\section{Quasiconjuagies}
  
\indent Suppose $F:M \rightarrow M$ and $f:N \rightarrow N$
are two dynamical systems. Recall that $F$ and $f$
are semi-conjugate if there is a map $\lambda:M \rightarrow N$
for which $f\circ \lambda= \lambda \circ F.$
 
\medskip

 It might happen that $dim(M)>dim(N$) and that $F$
mathematically models a situation more precisely than
$f.$ In that case we introduce the notion of 
$\underline{quasiconjugacy},$  $\Lambda: M \rightarrow N.$
For that, we decompose $M$ into three disjoint subsets, depending
on the (assumed) size of the difference between $f \circ \Lambda$
and $\Lambda \circ F.$
So we write, inspired by physics, $M=M_Q \cup M_{SC} \cup M_C$, and call
these
the $\underline{Quantum},\; \underline{Semi}-\underline{Classical}\;
 {\mbox{and}} \;
\underline{ Classical}$ regions.
We assume that the differences are largest in the Quantum region and
smallest in the Classical region. The system $(N,f)$ is called
${\underline{Classical}}$ and $(M,F)$ is said to be
the ${\underline{Quantization}}$ of $(N,f).$

\medskip

We say that $\Lambda$ satisfies the $\underline{Correspondence}\;
\underline{Principle}$ if the difference between $f^n \circ \Lambda$ and
$\Lambda \circ F^n$ is within some small bound (to be decided)
on $M_C$ while $f^n \circ \Lambda$ and $\Lambda \circ F^n$ have
sufficiently
close qualitative features on $M_{SC}$ 
(again in some sense to be made precise).
Here we might also put bounds on the number of iterates.

\medskip

In discussing quasiconjugacies, one often should think of
quasiconjugacy as a working hypothesis. In fact, the
conclusion of a research might rather be that a supposed quasiconjugacy
actually is not a quasiconjugacy or at least fails to be for
certain parameters.

\medskip

Let $u:M \rightarrow \RR$ (or more generally $\RR^n$)
and $v:N \rightarrow \RR$ be given functions, we call
 them $\underline{observables}$. The underlying working
hypothesis is then usually that $u$ and $v$ measure
the same phenomenon on the respective manifolds. We say that
$F$ $\underline{localizes}$ (with respect to $u,$ $v$) if for every
$x \in M$  the sequence $\{u(F^n(x))\}$ is bounded, while
there is a $y \in N$ for which
$\{v(f^n(y))\}$ 
is unbounded. This phenomenon has been experimentally
and numerically observed in some quantum systems $M$
for which the classical system $N$ is chaotic according
to numerical experiments. 
This then is considered a failure of the Correspondence Principle.
Another problem  for the correspondence principle has been found
in $\underline{resonances}$, where $u$ can be rigorously proved
to go to infinity at a faster pace than indicated numerically
for $v.$

\medskip

Here we discuss briefly an example arising from quantum mechanics
where $N$ is of complex dimension $1$
and $M$ is an infinite-dimensional complex manifold. Here $f$ is a
piecewise
holomorhic map
and $F$ is holomorphic. 

\medskip

For other rigorous works on related systems,
see Combescure ([Co]).

\medskip

Our example has the advantage that it is sufficiently simple that the
classical and quantized versions can be analyzed rigorously with
elementary means, as opposed to
the more complicated standard map ([Ch]) and its quantized counterpart,
the
kicked rotor ([IS]). Our main result is a mathematically rigorous proof
that
localization takes place. We believe this is the first such rigorous
proof,
previous results have been based on computer experiments.
Our main results is for the case which is considered furthest
from resonance, namely, with Golden Mean spacing between kicks.
Of course, it would be very interesting to extend this proof
to the kicked rotor case.

\section{The two dynamical systems}

\subsection{The $1-$dim complex dynamical system $(N,f)$}

\indent We will define a piecewise holomorphic map $f: N:=\CC^*
\rightarrow \CC^*,$ the punctured complex plane, $\CC^*
=\CC \setminus (0)$. Here $f$
will be the composition of two maps

$$
f=\Phi \circ \kappa
$$

\noindent where $\Phi$ is a rotation and $\kappa$ is a piecewise defined
scaling. 

\medskip

We set $\Phi(z)=e^{i \lambda} z$ for some real parameter $\lambda$ and

\bea
\kappa(z)  & = &  e^{-K}z \; \mbox{if Im}(z) \geq 0,\\
\kappa(z)  & = & e^K z \; \mbox{if Im}(z) <0\\
\eea

\noindent where $K$ is a nonzero real parameter. The set of points
with bounded orbit, i.e. bounded away from $0$ and $\infty$,
depends on $\lambda$ but not on $K.$

\medskip

Before introducing $M$ and $F$, we rewrite $f$ using variables
$\theta= {\mbox{arg}}\; z$ and $P=\log |z|$ and consider $\Phi$
and $\kappa$ as maps given by Hamiltonian flows on the cylinder
$T^1 \times \RR$ also called $N.$ This allows us to set up the
corresponding
Schrodinger equations.

\subsection{$(N,f)$ as a Hamiltonian system on a cylinder}

\indent Physically 
speaking, we model a particle moving on a circular, concentric
orbit in an  evenly charged spherical
plasma and being kicked by an external field
at given evenly spaced times.

\medskip

We use the variables $\theta$ and $P$ to denote angle
$\theta\in [0,2\pi)$ and angular momentum $P \in \RR$
respectively of the particle. Our manifold $N$ is $T^1 \times \RR.$
The particle moves with constant speed between kicks
with motion governed by the Hamiltonian $C(\theta,P)=\omega P$
where $\omega$ is the angular velocity. So the equations
of motion are

\bea
\frac{d\theta}{dt} & = & \frac{\partial C}{\partial P}\\
& = & \omega\\
\frac{d P}{dt} & = & -\frac{\partial C}{\partial \theta}\\
& = & 0\\
\eea

If the time between kicks is $T,$ this gives rise to a time-$T$
map on $N,$

$$
\Phi(\theta,P) = (\theta+\omega T,P).
$$ 

This corresponds to setting $\lambda=\omega T$.

\medskip

 We assume that the motion of the particle is governed by a Hamiltonian
$C'(\theta,P)=\frac{K H(\theta)}{\epsilon}$ during kicks, where
$H:T^1 \times \RR$ is continuous and piecewise smooth, $K$
is a parameter called the kick strength and $\epsilon$ is the duration of
the
kick. This gives rise to a time-$\epsilon$ map $\kappa: N \rightarrow N$
given by $\kappa(\theta,P)=(\theta,P-K H'(\theta)).$ Our example
corresponds
then to the choice of $H(\theta)$ being the tent map, normalized so
that $\int_0^{2\pi} H(\theta)=0.$

\bea
 H(\theta) & = &  \theta-\frac{\pi}{2}\; \mbox{if}\;0 \leq
 \theta \leq \pi,\\
H(\theta) & = & \frac{3\pi}{2}-\theta \; \mbox{if} 
\pi \leq \theta \leq 2 \pi\\
\eea

\subsection{The $\infty-$dim complex dynamical system $(M,F)$}

\indent The complex manifold $M$ consists of projective Hilbert space
$\PP L^2([0,2\pi])$, where $L^2([0,2\pi])$ is the
complex Hilbert space of complex valued $L^2$ functions. We can identify
$M$ with the boundary of the unit ball in $L^2([0,2\pi])$ when
we identify $\psi$ and $e^{i \theta} \psi.$

\medskip

The map $\Phi:N \rightarrow N$ corresponds to a map
$\Phi':M \rightarrow M$ obtained by solving the Schrodinger
equation with Hamiltonian $C(\theta,P)=\omega P.$ Since
$P$ corresponds to the operator $-i\hbar \frac{\partial}{\partial \theta}$
on $L^2$,
we get for $L^2$ functions on $[0,2\pi):$

\bea
i \hbar \frac{\partial \psi}{\partial t} & = & -i\hbar \omega \frac{
\partial \psi}{\partial \theta},\; {\mbox{so}}\\
\psi(\theta,t) & = & \psi(\theta-\omega t)\\
((\Phi')(\psi))(\theta) & = & \psi(\theta-\omega T)\\
\eea

Similarly the map $\kappa:N \rightarrow N$ corresponds to the map
$\kappa':M \rightarrow M$ obtained by solving the Schrodinger equation
with Hamiltonian $C'=\frac{K H(\theta)}{\epsilon},$ so

\bea
i \hbar \frac{\partial \psi}{\partial t} & = & \frac{K
H(\theta)}{\epsilon}
\psi,\; {\mbox{hence}}\\
\psi(\theta,t) & = & e^{-\frac{i K H(\theta) t}{\epsilon}}\psi(\theta)\\
((\kappa')(\psi))(\theta) & := & \psi(\theta,\epsilon) \\
& = & e^{-\frac{i K H(\theta)}{\hbar}}\psi(\theta).\\
\eea

Finally we define $F: M \rightarrow M$ by $F:=\Psi' \circ \kappa'.$

\subsection{The Quasiconjugacy}

\indent It is convenient to write $\psi\in M$ as Fourier series,
$\psi=\sum_{n=-\infty}^{n=\infty} a_n e^{in \theta}.$ The
$\{e^{in \theta}\}$ are then eigenfunctions for the operator
corresponding to $P$, with eigenvalues $\{n \hbar\}.$ Since
by the uncertainty principle $\Delta \theta \Delta P \geq
\hbar,$ if $\psi$ has most of its $L^2$ norm
in an angle $\Delta \theta \leq \frac{1}{20}$ say, then
$\Delta P \geq 20 \hbar.$ We use this to define the Quantum region
to consist of those $\psi$ for which $\sum_{|n|\leq 20}|a_n|^2
>.9.$ At the other extreme, when the particle
is in an eigenstate when $|n|>200$ say,
the particle is easily ionized by an external field, and one defines
the Classical region to consist of those $\psi$
for which $\sum_{|n|>200}>.9.$ The remaining $\psi$
belong to the Semi-Classical region. These definitions
are of course somewhat arbitrary, and should be modified
as appropriate in various cases.

\medskip

If the particle is kicked with a change in angular momentum on the
order of $\hbar$, and the kicks are random, then one can
expect that the angular momentum grows like $\sqrt{m}$ where
$m$ is the number of iterates. This means that a particle
in the Quantum region can be expected to ionize, i.e. 
reach the Classical region after about $4\times 10^4$
iterations. This then is a reasonable upper bound on the number of
iterations one ought to consider.

\medskip

Next we define the semi-conjugacy $\Lambda: M \rightarrow N$.
So, given $\psi \in M$, we need to associate a $(\theta,P)=
(\theta(\psi),P(\psi))=\Lambda(\psi)\in N.$ This is not canonical,
for example, if $\psi \equiv \frac{1}{\sqrt{2\pi}},$ the mass
is evenly distributed on the circle. (There is no canonical
center on a compact manifold.) We define

$$
\theta(\psi)=\sup_{0<\zeta<2\pi}\{\int_0^{\zeta}|\psi|^2<\frac{1}{2}\},
$$

\noindent so the particle can be found with equal probability in
 $(0,\theta(\psi))$
and in \\
$(\theta(\psi),2\pi).$ Of course this selects $\theta=0$
as a privileged point, which is rather arbitrary. 
  
\medskip

It is easier to define $P(\psi).$ We set 

$$
P(\psi)=\sup_m \{\sum_{n \leq m} |a_n|^2 \leq \frac{1}{2}\}.
$$

We notice then that the map $\Lambda:M \rightarrow N$ is not onto
since the values of $P(\psi)$ always is an integer.

\medskip

Next we define the observables $u:M \rightarrow \RR$ and
$v:N \rightarrow \RR$ using the angular momentum:
For every $\psi\in M$, set $u(\psi)=
\sup_m \{\sum_{|n|\leq m}|a_n|^2 \leq \frac{1}{2}\}$ and
for $(\theta,P)\in N,$ set
$v(\theta,P)=|P|.$

\begin{remark}
In our example we will be able to study the long-term behaviour
of $u,v$ directly without a detailed analysis of $\Lambda.$
The example suffers from a common defect, namely
we let $n \rightarrow \infty,$ instead of restricting to physically
meaningful $n$ as indicated above. This aspect should be analyzed. Note
also that analysis is often done by letting $\hbar \rightarrow 0$. This is 
similarly physically meaningless since $\hbar$ is a constant. We keep
$\hbar$ fixed here.
\end{remark}

\subsection{The iterates.}

\indent We can write down formulas for the $n^{\mbox{th}}$ iterates.

\bea
H_n(\theta) & := & \sum_{k=1}^n H(\theta-k \omega T)\\
(F^n(\psi))(\theta) & = & e^{-\frac{i}{\hbar}K H_n(\theta)}\psi(\theta-
n \omega T)\\
\tilde{H}_n'(\theta) & := & \sum_{k=0}^{n-1}H'(\theta+k \omega T)\\
f^n(\theta,P) & = & (\theta+n\omega T,
P-K\tilde{H}_n'(\theta))
\eea

We will consider the case when $\frac{\omega T}{2\pi}=
\frac{\sqrt{5}+1}{2},$
the case of the Golden Mean rotation. This is the case where one might 
perhaps expect
the highest amount of diffusion in the angular momentum.

\section{The Main Theorems}

\subsection{Remarks on the Golden Mean}

\indent To carry out our investigation of the localization properties
we will at first recall some simple properties of the Golden Mean.
For ease of reference we include  proofs.
We write $r=\frac{\sqrt{5}+1}{2}$ as the limit of the
fractions $r_n=\frac{p_n}{q_n}$ where
$\{q_n\}_{n=1}^\infty={2,3,5,8,\dots}$
is a Fibonacci sequence and $p_n=q_{n+1}.$
Note that $q_{3n+1}$ is even, the others are odd
and $r_{n+1}-r_n= \frac{(-1)^{n+1}}{q_{n+1}q_n}.$
Notice that for large $n$ we have the estimate
$|r_n-r|<c/(q_n)^2$ for some constant $c<1/2.$ 
We have that $r_n<r<r_k$ for $n$ even and $k$ odd.

Consider the map on the unit circle, $T$: $f(\theta)=\theta+2 \pi r.$

\begin{lem}
Let $\theta \in T.$ Then for each large $n$ the number
of iterates $\{f^k(\theta)\}_{k=1}^{q_n}$ in the upper
half circle $[0,\pi]$ (or $(0,\pi)$) is in the interval $[\frac{q_n}{2}-3,
\frac{q_n}{2}+3].$
\end{lem}

{\bf Proof:}
Divide the interval $[0,2\pi]$ into $q_n$ equal intervals\\
$I_j=[\frac{2\pi j}{q_n},  \frac{2\pi (j+1)}{q_n}),\; j=0,\dots q_n-1.$
Notice that if we ignore the error between $r_n$ and $r$,
the map simply is a permutation of the $I_j.$
If there was no error, then the number of points
in the orbit in the upper half circle is exactly half,
except that $q_n$ is sometimes odd and some points land exactly on
the endpoints of the intervals. However, the errors between $r$ and
$r_n$ are very small and even after $q_n$ iterates
they are less than half the length of the intervals,
i.e. at most $2\pi c/q_n.$ This makes at most an error
of $\pm 1$ in the counting.\\

\CQFD\\
   
\begin{lem}
Any positive integer $k$ can be written as a strictly
increasing sum $q_{j_1}+\dots + q_{j_\ell}$ of Fibonacci numbers
where $\ell \leq   \frac{\ln k}{\ln 1.5}.$
\end{lem}

{\bf Proof:} We prove this estimate by induction. Note that
$1.5 q_n \leq q_{n+1} \leq 2 q_n.$
Let $q_\ell$ be the largest Fibonacci number less than or equal
to $k.$ 
Then 

\bea
k-q_\ell  & \leq &  q_{\ell+1}-q_\ell\\
& = & q_{\ell-1}\\
& \leq &  \frac{1}{1.5} q_\ell\\
&  \leq &  \frac{1}{1.5} k.
\eea

\CQFD\\

Hence we can immediately extend the above estimate on the number of
points on an orbit in the upper half circle:

\begin{lem}
Let $k$ be any positive integer. Then for any $\theta\in [0,2\pi]$
the number of points in $[0,\pi]$ (or $(0,\pi)$) of the orbit
$\{f^n(\theta)\}_{1\leq n \leq k}$ differs from $k/2$ by
at most $3 \frac{\ln k}{\ln 1.5}.$
\end{lem} 

Our next goal is to estimate the functions $H_{q_n}$ and their
derivatives. 

\begin{lem}
$H_{q_n}$ is a Lip $1$ function, $|H'_{q_n}| \leq 3$ and
$|H_{q_n}(\theta)| \leq C \frac{n}{1.5^n}$
for some fixed constant $C.$
\end{lem}

{\bf Proof:} We carry out the details in the case $r_n<r$ and $q_n$ even,
the other cases are similar. For this, observe that at least for
large $n$ there exist $0<\alpha_j,\beta_j< \frac{1}{2q_n^2}$ and for
each $j=1,\dots,q_n$ corresponding $0 \leq \ell=\ell(j),
k=k(j)\leq q_n-1$ so that

\bea
f^\ell([\frac{\pi (j+1)}{q_n}-\alpha_j, \frac{\pi (j+1)}{q_n}])
& = & [\pi,\pi+\alpha_j]\\
f^k([\frac{\pi (j+1)}{q_n}-\beta_j, \frac{\pi (j+1)}{q_n}])
& = & [0,\beta_j]\\
\eea

Next we calculate the difference between $\ell$ and $k$.
But note that since $r=\frac{p_n}{q_n}+\delta$ with very small
$\delta$, we can write $\frac{q_n}{2}r= \frac{p_n}{2}+ \delta'$
where $\frac{p_n}{2}=\frac{1}{2} \;{\mbox{mod}}\; 1$ since necessarily
$p_n$ is odd. Hence we know that $|\ell-k|=\frac{q_n}{2}.$
In particular this means that $\alpha_j-\beta_j = \pm \frac{q_n}{2}
(r_n-r).$
Also, there must necessarily be equally many positive as negative ones
since on the average equally many points are in the upper as lower half
circle.

\medskip

If $\alpha_j>\beta_j,$ it means that $\ell=k+\frac{q_n}{2}$ and that
on the interval $\hat{I}_j:=\left(\frac{\pi(j+1)}{q_n}-\alpha_j,
 \frac{\pi(j+1)}{q_n}-\beta_j \right),$ the function $H'_{q_n}=-2.$ We say
that $\hat{I}_j$ is negative.   
If $\alpha_j<\beta_j,$ it means that $k=\ell+\frac{q_n}{2}$ and that
on the interval $\hat{I}_j:=\left(\frac{\pi(j+1)}{q_n}-\beta_j,
 \frac{\pi(j+1)}{q_n}-\alpha_j \right),$ the function $H'_{q_n}=2.$ We say
that $\hat{I}_j$ is positive.

\medskip

We want to prove the estimate of $H_{q_n}$. Since $\int |H'_{q_n}| \sim
1,$
which is too large for our purposes, we need to show
 that there is a lot of cancellations in order to
prove that $\max_\theta |\int_0^\theta H'_{q_n}|$ is small.
For this, define a string $S(i,j),0 \leq i<j<q_n$
to be the collection of intervals $\hat{I}_i,\dots,\hat{I}_j.$

\medskip

We want to show that for any string there is about
the same number of positive and negative intervals. For this we need
to understand the behaviour under iteration of nearest
neighbors in a string.

\medskip

Note that $r_{n-1}>r.$
Consider the iterate 

$$
f^{q_{n-1}}([\pi,\pi+\frac{2\pi}{q_n}])
=[\pi-2\pi q_{n-1}(r_{n-1}-r), \pi+\frac{2\pi}{q_n}-2\pi
q_{n-1}(r_{n-1}-r)].
$$

Observe that  
$\pi-2\pi q_{n-1}(r_{n-1}-r)=\pi-\frac{2\pi}{q_n}+2\pi q_{n-1}(r-r_n)$.
This implies that $f^{q_{n-1}}$ maps $I_j$ to $I_{j-1}$ with the
"usual" error coming from $r_n\neq r.$ Note that
$\hat{I}_j$ is negative when $\ell(j)\in [\frac{q_n}{2}+1,\dots,q_n]$
and $\hat{I}_j$ is positive when $\ell(j)\in [1,\dots,\frac{q_n}{2}].$

Hence $\hat{I}_j$ is negative if $\frac{2\pi}{q_n}\ell(j)\in (\pi,2\pi]$
and $\hat{I}_j$ is positive if  $\frac{2\pi}{q_n}\ell(j)\in (0,\pi]$.

\medskip

Next we compare  $\frac{2\pi}{q_n}\ell(j)$ and
$\frac{2\pi}{q_n}\ell(j-1),0 \leq j <q_n.$

\bea
f^{\ell(j)}\left(\frac{\pi(j+1)}{q_n}\right) & = & \pi+\alpha_j\\
f^{\ell(j-1)}\left(\frac{\pi j}{q_n}\right) & = & \pi+\alpha_{j-1}.\\
\eea

Also $f^{q_{n-1}}$ maps $I_j$ to $I_{j-1}.$ Hence as points on the unit 
circle,

\bea
\frac{2\pi [\ell(j-1)+q_{n-1}]}{q_n}  & = &  \frac{2\pi}{q_n}\ell(j)\\
\frac{2\pi \ell(j-1)}{q_n}+2\pi r_n  & = &  \frac{2\pi}{q_n}\ell(j).\\
\eea

Therefore, given a string $S(i,j)$, we can count positive and negative
intervals as numbers of points on the orbit of $\frac{2\pi \ell(i)}{q_n}$
in the upper and lower half circle of rotation by the 
Golden Mean. Hence Lemma 4.3 implies that
the difference between the number of positive and negative intervals
in a string $S(i,j)$ is at most $\frac{3 \ln q_n}{\ln 1.5}.$
Therefore

\bea
\left| \int_{S(i,j)} H_{q_n}\right| & \leq & \frac{3 \ln q_n}{1.5}
\frac{1}{2 q_n}\\
\leq C \frac{n}{(1.5)^n}\\
\eea

\CQFD\\

\subsection{Localization in the Quantized Case}

\indent Now we are ready to prove that the quantized map has an empty
basin of attraction at infinity in the topology
defined by $u$, i.e. that the map
localizes:

\begin{theo}
For any $\psi \in L^2([0,1]), \; \|\psi\|=1,$ and any
$\epsilon >0$ there exists an integer $N$ so that for any
iterate $F^n(\psi)=\sum_m c^n_m e^{i m \theta}$ we have
$\sum_{|m|<N} |c^n_m|^2 >1-\epsilon.$ In particular,
$\{u(F^n(\psi))\}_n$ is a bounded sequence for any $\psi.$
\end{theo}

{\bf Proof:}
Pick any $\delta>0$. Write $\psi=\sum_m c_me^{im\theta}.$
Let $N_\psi$ be chosen so that $\sum_{|m|<N_\psi} |c_m|^2 >1-\delta.$
Suppose that $h(\theta)$ is any continuous real function on $T^1$, and
suppose $\lambda=\{\lambda_m\}$ is any sequence of real numbers. We
denote by $\psi_\lambda$ the function $\sum_m c_m e^{i \lambda_m} e^{im
\theta}
$. 
Set $\tilde{\psi}=e^{ih}\psi_\lambda=\sum_m \tilde{c}_m e^{im\theta}.$
Then

\bea
\tilde{c}_m & = & \frac{1}{2\pi} \int_0^{2\pi}\tilde{\psi}
e^{-im \theta}\\
           & = &  \frac{1}{2\pi} \int_0^{2\pi} e^{ih}\psi_\lambda
e^{-im \theta}\\
 & = & c_m e^{i \lambda_m}+\frac{1}{2\pi} \int_0^{2\pi}
(e^{ih}-1)\psi_\lambda
e^{-im \theta},\\
\eea

\noindent so 

\bea
|\tilde{c}_m-c_m e^{i \lambda_m|} & \leq &
 \sup |h| \|\psi_\lambda\|_{L^1} \\
& \leq &  \sqrt{2\pi}  \sup |h| \|\psi_\lambda\|\\
& =  &  \sqrt{2\pi}  \sup |h| \|\psi\|\\
& =  &  \sqrt{2\pi}  \sup |h| \\
\eea

Hence if $\sup |h|<\eta$ for some constant
$\eta=\eta(\delta)>0$ 

$$
\sum_{|m|<N_\psi}|\tilde{c}_m|^2 \geq 1-2\delta.
$$

\medskip

There is an integer $k=k(\delta)$ so
 that $\frac{K}{\hbar}\sum_{\ell=k}^\infty \sup |H_{q_\ell}|
<\eta.$

\medskip

Next let $n$ be any integer. By Lemma 4.2 we can write

\bea
n & = & q_{n_1}+\cdots +q_{n_j}+q_{n_{j+1}}+ \cdots+q_{n_i}\\
& = &  q_{n_1}+\cdots +q_{n_j}+n'\\
& = & n''+n'\\
\eea

\noindent where $q_{n_1}<q_{n_2} < \cdots$ are Fibonacci numbers
and $n_j<k,n_{j+1} \geq k.$ We can then write

\bea
F^n(\psi) & = & F^{n''}(F^{n'})(\psi))\\
& = & F^{n''}(\tilde{\psi})\; {\mbox{if we define}}\\
\tilde{\psi} & = & F^{n'}(\psi). \; {\mbox{From Section 3.5:}}\\
\tilde{\psi} & = & e^{-\frac{i}{\hbar}H_{n'}(\theta)}
\psi(\theta-n' \omega T)\\
& = & e^{-\frac{i}{\hbar}
\sum_{\ell=j+1}^i H_{q_\ell}(\theta-\sum_{r=j+1}^{\ell-1}q_{n_r}
\omega T)}
\sum_m c_m e^{-i q_{n'}\omega T}e^{im \theta}\\
\eea

Setting  

\bea
h & = & -\frac{K}{\hbar} \sum_{\ell=j+1}^i H_{q_\ell}(\theta
-\sum_{r=j+1}^{\ell-1} q_{n_r} \omega T )\\
\lambda =\{\lambda_m\} & = & \{-m n' \omega T\}\\
\tilde{\psi} & = & e^{ih} \psi_\lambda \\
& = & \sum \tilde{c}_m e^{im \theta}\; {\mbox{we get:}}\\
\sum_{|m|<N_\psi} |\tilde{c}_m|^2 & \geq & 1-2\delta.\\
\eea

\medskip

Note that the collection of all possible operators $F^{q_{n_1}+\cdots
q_{n_j}}$
is finite.

\medskip

We write $\tilde{\psi}=\tilde{\psi}_1+\tilde{\psi}_2$
where $\tilde{\psi}_1=\sum_{|m| < N_\psi}\tilde{c}_m
e^{im\theta}.$ Then $\|\tilde{\psi}_2\| \leq \sqrt{2\delta}$ and
$\|F^{q_{n_1}+\cdots +q_{n_j}}(\tilde{\psi}_2)\|\leq \sqrt{2\delta}.$

Also, by compactness of the space $\{\{s_m\}_{|m|\leq N_\psi};
1-2\delta \leq \sum|s_m|^2 \leq 1\},$ it follows that for some
$N$, if $(F^{q_{n_1}+\cdots +q_{n_j}})(\tilde{\psi}_1)=\sum
\hat{c}_m e^{im\theta}$ then $\sum_{|m|< N}|\hat{c}_m|^2 \geq 1-3\delta.$

\medskip

It follows that 

$$
\left(\sum_{|m|>N}|\hat{c}_m|^2 \right)^{\frac{1}{2}}<
\sqrt{3\delta},
$$

\noindent so

$$
\left(\sum_{|m|>N}|c_m^n|^2 \right)^{\frac{1}{2}}<
\sqrt{3\delta}+\sqrt{2\delta}
$$

\noindent 

\noindent and hence

\noindent $\sum_{|m|< N} |c_m^n|^2 \geq 1-12\delta.$ Set
$\delta=\epsilon/12$
to complete the proof.\\

\CQFD\\

\subsection{Diffusion in the Classical Case}
 
\indent Next we show that the classical dynamics does not localize. We
show
(slow) diffusion to infinity, which is stronger than
finding one unbounded orbit. Diffusion is sometimes considered as
Chaos in the Classical System, while the Localization
in the previous subsection is considered as a Supression of Chaos
in the Quantized Case, assumed to be a result of Destructive
Interference.

\begin{theo}
Let $I=\{(0,2\pi)\times (0)\}$ be a set of initial states
for the classical system. Let $\epsilon>0$ and $N>0$ be given.
Then there exists an integer $n$ so that $\{\theta; f^n(\theta,0)
\in \{|P|<N\}\}$ has linear measure at most $\epsilon.$
\end{theo}

{\bf Proof:}
We will first discuss a situation where we consider that
$\frac{|\alpha_j-\beta_j|}{\frac{2\pi}{q_n}}=\delta$ is independent of
$n$ and then we correct for the error at the end.

\medskip

So suppose at first that $0<\delta <1/2$ is given. Let
$\chi_n(\theta)$ be an integer-valued step function with
values in $\{-n,\dots,n\}.$ Set $E_{n,k}=\{\theta;
\chi_n=k\}.$ We assume that ${\chi_{n+1}}_{|E_{n,k}}$ is a stepfunction 
with values in $\{k-1,k,k+1\}.$ The set of points
with value $k$ in $E_{n,k}$ has measure $(1-2\delta)|E_{n,k}|$
and the set of points with value $k\pm 1$ both have
measure $\delta |E_{n,k}|.$ We get inductively
that $|E_{n,k}|=\sum_{i-j=k}a_{i,j}$ where
the coefficients $a_{i,j}$ are given by the binomial
formula:

$$
((1-2\delta)+\delta x+\delta y)^n=\sum a_{i,j}x^i y^j.
$$
 
In this simplified case, the theorem follows from a simple
calculation which shows:

\begin{lem}
For any given $k,$ $\lim_{n\rightarrow \infty}|E_{n,k}|=0.$
\end{lem}

Finally to take into account that $\epsilon$ actually will vary,
observe that 

$$
\frac{|\alpha_j-\beta_j|}{\frac{2\pi}{q_n}}=\frac{q_n^2|r_n-r|}{2\pi}
=:\delta_n
$$

\noindent stays bounded and bounded away from zero. We may assume by
taking
a subsequence that $\delta_n \rightarrow \delta$ arbitrarily fast.
Pick $m_0$ large enough that $|E_{n,k}|<\frac{\epsilon}{3N}$ whenever
$k<N$ and $m \geq m_0.$ Then we choose Fibonacci numbers
$1<<q_{k_1} << q_{k_2}<< \cdots << q_{k_{m_0}}$ for which
$r_{k_j}<r$ and the $q_k$ are even (as we discussed in the proof
of Lemma 4.4). Then $n=\sum q_{k_j}$ will work.\\

\CQFD\\

\begin{coro}
Almost every sequence $\{v(f^n(\theta,P))\}$ is unbounded. In particular,
there does exist a $y\in N$ for which all $f^n(y)$ are well defined
and $\{v(f^n(y))\}$ is unbounded.
\end{coro}

%%%%%%%%%%%%%%%%%%%%%%%%%%%%%%%%%%%%%%%%%%%%%%%%
%%%%%%%%%%%%%%%%%%%%%%%%%%%%%%%%%%%%%%%%%%%%%%%%

\noindent John Erik Fornaess\\
Mathematics Department\\
The University of Michigan\\
East Hall, Ann Arbor, Mi 48109\\
USA\\

\end{document}